\documentclass[reqno]{amsart}
\usepackage{amsmath,amssymb,colordvi}
\usepackage{multicol}
\usepackage{color,bbm}
\usepackage{epsfig}

 \newtheorem{theorem}{Theorem}[section]
 
 \newtheorem{lemma}[theorem]{Lemma}

 \newtheorem{definition}[theorem]{Definition}
 \newtheorem{remark}[theorem]{Remark}
 
 \numberwithin{equation}{section}

 \usepackage[dvips,bottom=1.4in,right=1in,top=1in, left=1in]{geometry}

\begin{document}

\title[Exponential Stability of Fractional Stochastic System ]{Exponential Stability of Nonlinear Fractional Stochastic System with Poisson Jumps}

\author{P. Balasubramaniam}
\address{%
P. Balasubramaniam, Department of Mathematics,
The Gandhigram Rural Institute - Deemed to be University,
Gandhigram - 624 302, Dindigul, Tamil Nadu, India.}
\email{balugru@gmail.com}

\author{T. Sathiyaraj}
\address{%
T. Sathiyaraj, Institute of Computer Science \& Digital Innovation,  UCSI University, Malaysia.}
 \email{sathiyaraj133@gmail.com }

\author{K. Priya}
\address{ K. Priya, Department of Mathematics,
	The Gandhigram Rural Institute - Deemed to be University,
	Gandhigram - 624 302, Dindigul, Tamil Nadu, India. }
\email{priyak250796@gmail.com}

\thanks{University Grants Commission - Special Assistance Program (Department of Special Assitance -I), New Delhi, India, File No. F. 510/7/DSA-1/2015(SAP-I)}



\begin{abstract}
In this paper  exponential stability of nonlinear fractional order stochastic system with Poisson jumps is studied in finite dimensional space. Existence and uniqueness of solution, stability and exponential stability results are established by using boundedness properties of Mittag-Leffler matrix function, fixed point route and local assumptions on nonlinear terms. A numerical example is given to illustrate the efficiency of the obtained results. Finally, conclusion is drawn.
\end{abstract}

\maketitle

\section{Introduction}
Fractional differential equations (FDEs) have attracted much attention and have been widely used in engineering, physics, chemistry, biology, and other fields \cite{1}. The FDEs  have captured many researchers \cite{5, 2,4,3} through its potential applications.
The concept of stability is most important qualitative property of dynamical control systems since every control system must be primarily stable and only then other properties can be studied. Unstable systems have no practical significance in the real world situations. The key problem is to obtain information on the system's behavior (solution trajectory) without solving the differential equations. The theory considers the system’s behavior over a long period of time, that is, as $t\rightarrow\infty.$ The study of stability analysis of fractional order systems can be characterized by its behavior of the solutions. The question of stability is an interesting one and its study is meaningful in fractional order stochastic systems. In 1992, the general concept of stability has been defined by  Lyapunov \cite{24} for both linear and nonlinear systems. 

Initially, Matignon \cite{R2} presented the qualitative properties of stability and asymptotic stability results for the following autonomous system
$$\frac{d^\alpha x}{dt^\alpha}=Ax,\quad x(0)=x_0.$$
He guaranteed stability results by roots of the polynomials lie outside of the closed angular sector $|arg(\lambda)|\leq\frac{\alpha \pi}{2}.$ Further,  the above linear autonomous system has been extended to the nonlinear system (see \cite{R1}).

The qualitative behavior of stability and asymptotic stability results have been obtained by using Lipschitz condition.  Moreover,
by utilizing the fractional order comparison principle, authors in \cite{6,7} studied the stability of nonlinear FDEs. In \cite{8}, authors introduced the Mittag-Leffler stability of nonlinear FDEs with order $\alpha\in(0,1).$ The fixed point theorems have been used to study the stability of integer order differential system by many authors, notably Burton \cite{10}. The stability of Caputo type FDEs of order $\alpha \ (0<\alpha<1)$ has been considered by Burton and Zhang in \cite{11}, via resolvent theory and fixed point theorems in a weighted Banach space.

The concept of exponential stability plays a crucial role in the dynamical system and its convergence rate is faster than the asymptotic stability.
Further, existence and uniqueness of solutions for stochastic system with Poisson jumps \cite{R6}, exponential stability for stochastic partial differential equation \cite{R7}, exponential stability for stochastic system with semimartingales \cite{RR8}, existence and stability results for second order stochastic system with fractional Brownian motion \cite{R5}, existence, uniqueness of mild solution and stability results for second order stochastic system with delay and Poisson jumps \cite{R3},
exponential stability for second order stochastic system with Poisson jump \cite{R4} and exponential stability for stochastic system driven by G-Brownian motion \cite{R8} have emerged in the literature.
In \cite{64} existence and exponential stability of solution of nonlinear stochastic partial integrodifferential equations have been studied.  

However, it should be further emphasized that the stability and exponential stability of fractional stochastic system with Poisson jumps in finite dimensional space is yet to be elaborated, compared to that of integer order stochastic system. In this regard, it is necessary and important to study the exponential stability of nonlinear fractional stochastic system with Poisson jumps. To the best of   authors knowledge, only a few steps are taken to use the fixed point theorems to investigate the stability of fractional stochastic system. There is no work reported in the literature to study the problem of stability and exponential stability of fractional stochastic system with Poisson jumps in finite dimensional space. In order to fill-up this gap, in this paper, we study the following nonlinear fractional stochastic system with Poisson jumps
\begin{align}\label{p}
^CD_{0+}^q x(t)=& Ax(t)+f(t,x(t))+\int _0^t\sigma(s,x(s))dw(s)+\int _{-\infty}^{+\infty} g(t,x(t),\eta)\tilde{N}(dt,d\eta),\ t\in J\nonumber \\ 
x(0)=& x_0 
\end{align}
where $J:=[0,T],\ \frac{1}{2}<q<1$, $A$ is an diagonal stability matrix, nonlinear functions\ $f:J\times \mathbb{R}^n\rightarrow \mathbb{R}^n$,~ $\sigma:J \times \mathbb{R}^n\rightarrow \mathbb{R}^{n\times n}$ and $g:J \times \mathbb{R}^n \times \mathbb{R}_+ \rightarrow \mathbb{R}^n$ are continuous. 
Let $\{N(t),t\in J\}$ is a Poisson point process taking its values in a measurable space with a characteristic measure of $\varPi.$ Denote $N(dt,d\eta)$ as a Poisson counting measure induced by a Poisson pointing process $N(\cdot).$  $\tilde{N}(dt,d\eta)=N(dt,d\eta)-\varPi(d\eta)dt$ stands for the compensating measure  independent of the $n-$dimensional Wiener process ${w(t)}$ and satisfies $\int_{-\infty}^{+\infty}\varPi(d\eta)<\infty.$ 
This section ends by highlighting the  main contributions of this paper:
\begin{itemize}
	\item Existence and uniqueness of solution for nonlinear fractional stochastic system is proved.
	\item Stability results are investigated in finite dimensional stochastic settings.
	\item Exponential stability results are established using appropriate hypotheses.
	\item Numerical simulation is given to validate the proposed theoretical results.
\end{itemize}	

The outline of this paper is described as follows: In Section 2, preliminaries are given. In Section 3, stability and exponential stability results are established for nonlinear fractional order stochastic system with Poisson jumps by using the  Banach contraction mapping and Gronwall's inequality. We give an example to illustrate the theoretical results in Section 4. Finally, conclusion is drawn in Section 5.
\section{Preliminaries}
Let $(\Omega,\mathcal{F},P)$ be the complete probability space with a probability measure $P$ on $\Omega$ and
let $\left\lbrace \mathcal{F}_t|t\in J\right\rbrace $ be the filtration generated by $\left\lbrace w(s): 0\leq s\leq t\right\rbrace $ defined on the probability space $(\Omega,\mathcal{F},P)$. Let $L_2(\Omega,\mathcal{F}_t,\mathbb{R}^n)$ denotes the Hilbert space of all $\mathcal{F}_t$-measurable square integrable random variables with values in $\mathbb{R}^n$. 
Let $\mathcal{K}$=$\left\lbrace x(t):x(t)\in \mathcal{C}(J,L_2(\Omega,\mathcal{F}_t,\mathbb{R}^n))\right\rbrace$ be a Banach space  of all continuous square integrable and $\mathcal{F}_t$-adapted process with the norm $\|x\|^2=\sup\limits_ {t\in J}\mathbb{E}\|x(t)\|^2$, and
$\xi(\epsilon)=\{x:x \in \mathcal{K}:\mathbb{E}\|x\|^2\leq\epsilon,\ \epsilon >0\}.$ 
\begin{definition}{\cite{1}}
	Riemann-Liouville fractional operators:
	\begin{align*}
	(I_{0+}^qf)(x)=&\frac{1}{\Gamma (q)}\int\limits_0^x (x-t)^{q-1}f(t)dt \\
	(D_{0+}^qf)(x)=& D^n(I_{0+}^{n-q}f)(x),
	\end{align*}
	and the Laplace transform of the Riemann-Liouville fractional integral is given by 
	\begin{equation*}
	\mathcal{L}\{I_{t}^qf(t)\}=\frac{1}{\lambda^q}
	\hat{f}(\lambda),
	\end{equation*}
	where
	$\hat{f}(\lambda)= \int\limits_0^\infty e^{-\lambda t} f(t)dt.$
\end{definition}
\begin{definition} {\cite{1}}
	Caputo fractional derivative of order $n-1<q<n$ for a function $f:[0,\infty)\rightarrow\mathbb{R}$ is defined as
	\begin{equation*}
	^CD_t^q f(t)=\frac{1}{\Gamma(n-q)}\int\limits_{0}^t
	\frac {f^n(s)}{(t-s)^{q-n+1}}ds 
	\end{equation*}
	and its Laplace transform is
	\begin{equation*}
	\mathcal{L}\{^CD_t^q f(t)\}(s)= s^q f(s)-\sum_{k=0}^{n-1}f^k(0^+)s^{q-1-k}.
	\end{equation*}		
\end{definition}	
\begin{definition} {\cite{1}} Mittag-Leffler matrix function:\\
	A two parameter Mittag-Leffler matrix function type is defined by the following series expansion:
	\begin{align*}
	E_{q,p}(z)=\sum^\infty_{k=0}\frac{z^k}{\Gamma(kq+p)}, \quad \ q,p>0,\quad z\in\mathbb{C}.
	\end{align*}
	The most interesting properties of Mittag-Leffler function are associated with their Laplace integral
	\begin{align*}
	\int^\infty_0 e^{-st}t^{p-1}E_{q,p}(\pm at^q)dt = \frac{s^{q-p}}{(s^q\mp a)}.
	\end{align*}
	That is,\quad $\mathcal{L}\{t^{p-1}E_{q,p}(\pm at^q)\}(s)=\frac{s^{q-p}}{(s^q \mp a)}.$
\end{definition}
For more details on stochastic settings in finite dimensional space, Riemann-Liouville fractional derivative, Caputo fractional derivative and Mittag-Leffler function one can refer the papers \cite{70,sb2} and the references cited therein.
\begin{definition}
	Trivial solution $x(t)\equiv 0$ of (\ref{p}) is said to be stable in $\mathcal{K}$, if for any given  $\epsilon>0$ there exists $\delta(\epsilon)>0$ such that  $\mathbb{E}\|x_0\|\leq\delta$ satisfies $\mathbb{E}\|x(t)\|\leq \epsilon$ for all $t\geq 0.$ 
\end{definition}
\begin{definition}
	System $({\ref {p}})$ is said to be exponentially stable if  there exists two positive constants $\mu>0$ and $\mathbb{M}^* \geq 1$ such that
	\begin{align*}
	\mathbb{E} \|x(t)\|^2\leq \mathbb{M}^* e^{-\mu t}, \ \ \ t\geq 0. 
	\end{align*} 
\end{definition}
\begin{lemma}\cite{71}
	Let $p\geq 2$ and $g\in \mathcal{M}^2(J,\mathbb{R}^{n\times n})$ such that 
	\begin{align*}
	\mathbb{E}\int_0^T|g(s)| ^p ds <\infty.
	\end{align*}
	Then
	\begin{align*}
	\mathbb{E}\left| \int_0^T g(s)dB(s)\right|^p \leq\Big(\frac{p(p-1)}{2}\Big)^\frac{p}{2} T^\frac{p-2}{2} \mathbb{E}\int\limits_0^T|g(s)|^p  ds .
	\end{align*}
\end{lemma}
\begin{lemma}\cite{R20}
	For any $p\geq2$ there exists $\overline{c_p}>0$ such that
	\begin{align*}
	\mathbb{E}\sup\limits_{s\in [0,t]}\left\|\int^s_0\int^{+\infty}_{-\infty}g(\nu,z)\widehat{N}(d\nu,dz)\right\|^p&\leq  \overline{c_p}\Bigg\{\mathbb{E}\left[\left(\int^t_0\int^{+\infty}_{-\infty}\|g(s,z)\|^2\kappa (dz)ds\right)^\frac{p}{2}\right]
	\\&\quad+\mathbb{E}\left[\int^t_0\int^{+\infty}_{-\infty}\|g(s,z)\|^p\kappa(dz)ds\right]\Bigg\}.
	\end{align*}
\end{lemma}		
Solution of the system (\ref{p}) can be described as follows (see \cite{70})
\begin{align*}
x(t)=& E_q(At^q)x_0+\int_0^t(t-s)^{q-1}E_{q,q}(A(t-s)^q) \left[f(s,x(s))+\int _0^s\sigma(\theta,x(\theta))dw(\theta)\right]ds\\
&+\int_0^t(t-s)^{q-1}E_{q,q}(A(t-s)^q)\int _{-\infty}^{+\infty} g(s,x(s),\eta)\tilde{N}(ds,d\eta).
\end{align*}
Let us take the following hypotheses for further discussions:
\begin{itemize}
	\item[$(H_1)$]~Functions $f$, $\sigma$ and $g$ are continuous and there exists a constant $\beta>1$ and the functions $L_f(\cdot), L_\sigma(\cdot)$ and $L_g(\cdot)\in L^{\beta}(J,\mathbb{R}^+)$ such that
	\item[(i)]$\|f(t,x)-f(t,y)\|^2 \leq L_f(t)\|x-y\|^2$
	\item[(ii)]$\|\sigma(t,x)-\sigma(t,y)\|^2 \leq L_\sigma(t)\|x-y\|^2$
	\item[(iii)]$\int_{-\infty}^{+\infty}\|g(t,x,\eta)-g(t,y,\eta)\|^2\varPi(d\eta) \leq L_g(t)\|x-y\|^2.$
\end{itemize}
\begin{itemize}
	\item[$(H_2)$] Following properties hold
	for $t,s\geq 0$ 
	\item[(i)]$\|E_q(At^q)\|\leq N_1 e^{-\omega t}, \quad N_1\geq 1$ 
	\item[(ii)]$\|E_{q,q}(A(t-s)^q)\|\leq N_2 e^{-\omega(t-s)}, \quad N_2\geq 1.$
\end{itemize}
\begin{itemize}
	\item[$(H_3)$]There exist the constants $\widehat{V}_f,\widehat{V}_\sigma\ \text{and}\  \widehat{V}_g$ such that
	\item [(i)]  $\mathbb{E}\|f(t,x)\|^2\leq \widehat{V}_f\ (1+\mathbb{E}\|x\|^2)$
	\item [(ii)]  $\mathbb{E}\|\sigma(t,x)\|^2\leq \widehat{V}_\sigma \  (1+\mathbb{E}\|x\|^2)$
	\item [(iii)]  $\int\limits_{-\infty}^{+\infty}\mathbb{E}\|g(t,x,\eta)\|^2\Pi(d\eta)\leq \widehat{V}_g\ (1+\mathbb{E}\|x\|^2).$
\end{itemize}
In addition, we set
\begin{align*}
1.\ Q_1=&8N_2^2\left[\frac{1-e^{-2\omega T\alpha }}{2\omega\alpha}\right]^\frac{1}{\alpha}
\Bigg(\frac{T^{2q-1}}{2q-1}\|L_f\|_{L^{\beta}(J,\mathbb{R}^+)}
+\frac{T^{2q}}{q^2}\|L_\sigma\|_{L^{\beta}(J,\mathbb{R}^+)}
+\frac{T^{2q-1}}{2q-1}\overline{c_p}\|L_g\|_{L^{\beta}(J,\mathbb{R}^+)}\Bigg)\\
2.\ Q_2=&8N_2^2\left[\frac{1-e^{-2\omega T}}{2\omega}\right]\Bigg(\frac{T^{2q-1}}{2q-1}R_f
+\frac{T^{2q}}{q^2}R_\sigma+\frac{T^{2q-1}}{2q-1}\overline{c_p}R_g\Bigg) 
\end{align*}
\section{Main results}
\begin{theorem}\label{th1} Suppose that hypotheses (H1) - (H2) are satisfied. Then, the  system (\ref{p}) has atleast one solution provided that
	\begin{align}
	M=& 3N_2^2\left[\frac{1-e^{-2\omega T\alpha }}{2\omega\alpha}\right]^\frac{1}{\alpha}\left(\frac{T^{2q-1}}{2q-1}\|L_f\|_{L^{\beta}(J,\mathbb{R}^+)}+\frac{T^{2q}}{q^2}\|L_\sigma\|_{L^{\beta}(J,\mathbb{R}^+)}+\frac{T^{2q-1}}{2q-1}\overline{c_p}\|L_g\|_{L^{\beta}(J,\mathbb{R}^+)}\right) <1\label{22}
	\end{align}
	where $\frac{1}{\alpha}+ \frac{1}{\beta}=1,\ \alpha,\beta>1$. The trivial solution $x\equiv 0$ of (\ref{p}) is stable in $\mathcal{K}.$  
	
	\noindent\textbf{Proof:}
	For each positive number $r$ define  $\mathcal{K}_r=\Big\{x \in \mathcal{K}:\mathbb{E}\|x\|^2\leq r\Big\}$. Then for each $r$, $\mathcal{K}_r$ is obviously a bounded, closed and convex subset of $\mathcal{K}$. Next, we set $R_f=\sup\limits_{t\in J} \mathbb{E}\|f(t,0)\|^2,R_\sigma=\sup\limits_{t\in J} \mathbb{E}\|\sigma(t,0)\|^2$ and $R_g=\sup\limits_{t\in J} \mathbb{E}\|g(t,0)\|^2.$
	Now, we define a operator $\Phi:\mathcal{K}_r\rightarrow\mathcal{K}_r$ as follows
	\begin{align*}
	(\Phi x)(t)=& E_q(At^q)x_0+\int_0^t(t-s)^{q-1}E_{q,q}(A(t-s)^q)\left[f(s,x(s))+\int _0^s\sigma(\theta,x(\theta))dw(\theta)\right]ds\\&
	+\int_0^t(t-s)^{q-1}E_{q,q}(A(t-s)^q)\int _{-\infty}^{+\infty} g(s,x(s),\eta)\tilde{N}(ds,d\eta).
	\end{align*}
	\textbf{Step: 1.}  We prove that there exists a positive number $r$ such that $\Phi(\mathcal{K}_r)\subseteq\mathcal{K}_r$.
	\begin{align}\label{s1}
	\mathbb{E}\|(\Phi x)(t)\|^2&\leq
	4\mathbb{E}\|E_q(At^q)x_0\|^2+4\mathbb{E}\left\|\int_0^t(t-s)^{q-1}E_{q,q}(A(t-s)^q)f(s,x(s))ds\right\|^2\nonumber\\
	&\quad+4\mathbb{E}\left\|\int_0^t(t-s)^{q-1}E_{q,q}(A(t-s)^q)\left[\int _0^s\sigma(\theta,x(\theta))dw(\theta)\right]ds\right\|^2\nonumber\\
	&\quad+4\mathbb{E}\left\|\int_0^t(t-s)^{q-1}E_{q,q}(A(t-s)^q)\int _{-\infty}^{+\infty} g(s,x(s),\eta)\tilde{N}(ds,d\eta)\right\|^2\nonumber\\
	&=I_1+I_2+I_3+I_4.
	\end{align}
	By using $(H_1)-(H_2)$ and H\"{o}lder inequality, we get
	\begin{align*}
	I_1&=4\mathbb{E}\|E_q(At^q)x_0\|^2
	\leq 4\| E_q(At^q)\|^2\mathbb{E}\|x_0\|^2\\
	I_1&\leq 4N_1^2 e^{-2\omega T}\mathbb{E}\|x_0\|^2\\
	I_2&=4\mathbb{E}\left\|\int_0^t(t-s)^{q-1}E_{q,q}(A(t-s)^q)f(s,x(s))ds\right\|^2\\
	&\leq4	\int_0^t\Big((t-s)^{q-1}\Big)^2ds\int_0^t\|E_{q,q}(A(t-s)^q)\|^2 \mathbb{E}\|f(s,x(s))\|^2ds\\
	&\leq4\frac{T^{2q-1}}{2q-1}N_2^2\int_0^te^{-2\omega(t-s)}\Big(	\mathbb{E}\|f(s,x(s))-f(s,0)+f(s,0)\|^2\Big)ds\\
	&\leq 8\frac{T^{2q-1}}{2q-1}N_2^2\Bigg(\int_0^te^{-2\omega(t-s)}L_f(s) ds \mathbb{E}\|x\|^2+\int_0^te^{-2\omega(t-s)}	\mathbb{E}\|f(s,0)\|^2ds\Bigg)\\
	&\leq 8\frac{T^{2q-1}}{2q-1}N_2^2
	\Bigg(r\left[\int_0^te^{-2\omega(t-s)\alpha}ds\right]^\frac{1}{\alpha}
	\left[\int_0^tL_f^{\beta}(s) ds\right]^\frac{1}{\beta}+R_f\int_0^te^{-2\omega(t-s)}ds\Bigg)\\
	I_2&\leq8\frac{T^{2q-1}}{2q-1}N_2^2
	\Bigg(r\left[\frac{1-e^{-2\omega T\alpha }}{2\omega\alpha}\right]^\frac{1}{\alpha}\|L_f\|_{L^{\beta}(J,\mathbb{R}^+)}+R_f\left[\frac{1-e^{-2\omega T}}{2\omega}\right]\Bigg).
	\end{align*} 
	Similarly,
	\begin{align*}
	I_3&=	4\mathbb{E}\left\|\int_0^t(t-s)^{q-1}E_{q,q}(A(t-s)^q)\left[\int _0^s\sigma(\theta,x(\theta))dw(\theta)\right]ds\right\|^2\\
	&\leq 8\frac{T^{2q}}{q^2}N_2^2
	\Bigg(r\left[\frac{1-e^{-2\omega T\alpha }}{2\omega\alpha}\right]^\frac{1}{\alpha}\|L_\sigma\|_{L^{\beta}(J,\mathbb{R}^+)}+R_\sigma\left[\frac{1-e^{-2\omega T}}{2\omega}\right]\Bigg)
	\end{align*}
	and
	\begin{align*}
	I_4&=4\mathbb{E}\left\|\int_0^t(t-s)^{q-1}E_{q,q}(A(t-s)^q)\int _{-\infty}^{+\infty} g(s,x(s),\eta)\tilde{N}(ds,d\eta)\right\|^2\\
	&\leq8\frac{T^{2q-1}}{2q-1}N_2^2\overline{c_p}
	\Bigg(r\left[\frac{1-e^{-2\omega T\alpha }}{2\omega\alpha}\right]^\frac{1}{\alpha}\|L_g\|_{L^{\beta}(J,\mathbb{R}^+)}+R_g\left[\frac{1-e^{-2\omega T}}{2\omega}\right]\Bigg).
	\end{align*}
	From the above estimations, equation (\ref{s1}) becomes 
	\begin{align*}
	&\mathbb{E}\|(\Phi x)(t)\|^2\\
	&\leq 4\Bigg\{N_1^2 e^{-2\omega T}\mathbb{E}\|x_0\|^2+2\frac{T^{2q-1}}{2q-1}N_2^2
	\Bigg(r\left[\frac{1-e^{-2\omega T\alpha }}{2\omega\alpha}\right]^\frac{1}{\alpha}\|L_f\|_{L^{\beta}(J,\mathbb{R}^+)}+R_f\left[\frac{1-e^{-2\omega T}}{2\omega}\right]\Bigg)\\
	&\quad+2\frac{T^{2q}}{q^2}N_2^2
	\Bigg(r\left[\frac{1-e^{-2\omega T\alpha }}{2\omega\alpha}\right]^\frac{1}{\alpha}\|L_\sigma\|_{L^{\beta}(J,\mathbb{R}^+)}+R_\sigma\left[\frac{1-e^{-2\omega T}}{2\omega}\right]\Bigg)\\
	&\quad+2\frac{T^{2q-1}}{2q-1}N_2^2\overline{c_p}
	\Bigg(r\left[\frac{1-e^{-2\omega T\alpha }}{2\omega\alpha}\right]^\frac{1}{\alpha}\|L_g\|_{L^{\beta}(J,\mathbb{R}^+)}+R_g\left[\frac{1-e^{-2\omega T}}{2\omega}\right]\Bigg)\Bigg\}\\
	&\leq 4N_1^2 e^{-2\omega T}\mathbb{E}\|x_0\|^2+8N_2^2r\left[\frac{1-e^{-2\omega T\alpha }}{2\omega\alpha}\right]^\frac{1}{\alpha}
	\Bigg(\frac{T^{2q-1}}{2q-1}\|L_f\|_{L^{\beta}(J,\mathbb{R}^+)}
	+\frac{T^{2q}}{q^2}\|L_\sigma\|_{L^{\beta}(J,\mathbb{R}^+)}\\
	&\quad+\frac{T^{2q-1}}{2q-1}\overline{c_p}\|L_g\|_{L^{\beta}(J,\mathbb{R}^+)}\Bigg)+8N_2^2\left[\frac{1-e^{-2\omega T}}{2\omega}\right]\Bigg(\frac{T^{2q-1}}{2q-1}R_f
	+\frac{T^{2q}}{q^2}R_\sigma+\frac{T^{2q-1}}{2q-1}\overline{c_p}R_g\Bigg)\\
	&\leq 4N^2_1e^{-2\omega T}\mathbb{E}\|x_0\|^2+Q_1r+Q_2=r
	\end{align*}
	for
	$$ r= \frac{4N^2_1e^{-2\omega T}\mathbb{E}\|x_0\|^2+Q_2}{1-Q_1}, \quad Q_1<1.$$
	
	Hence, we obtain $\Phi(\mathcal{K}_r)\subseteq\mathcal{K}_r.$
	
	\noindent\textbf{Step: 2.} We show that $\Phi$ is  a contraction mapping.
	
	Let $x,y\in\mathcal{K}_r.$ In viewing of $(H_1)$ and $(H_2),$ for each $t\in J,$ we have 
	\begin{align*}
	&\mathbb{E}\|(\Phi x)(t)-(\Phi y)(t)\|^2\\ \leq& 3 \mathbb{E}\left\|\int_0^t(t-s)^{q-1}E_{q,q}(A(t-s)^q)\Big[f(s,x(s))-f(s,y(s))\Big]ds\right\|^2\\
	&+3\mathbb{E}\left\|\int_0^t(t-s)^{q-1}E_{q,q}(A(t-s)^q)\left[\int_0^s\Big(\sigma(\theta,x(\theta))-\sigma(\theta,y(\theta))\Big) dw(\theta)\right]ds\right\|^2\\
	&+3\mathbb{E}\left\|\int_0^t(t-s)^{q-1}E_{q,q}(A(t-s)^q)\int_{-\infty}^{+\infty}\Big( g(s,x(s),\eta)- g(s,y(s),\eta)\Big)\tilde{N}(ds,d\eta)\right\|^2\\
	\leq& 3\frac{T^{2q-1}}{2q-1}N^2_2
	\left(\int_0^te^{-2\omega(t-s)\alpha}ds\right)^\frac{1}{\alpha}
	\left(\int_0^tL_f^{\beta}(s)ds\right)^\frac{1}{\beta}\mathbb{E}\|x(t)-y(t)\|^2\\
	&+3\frac{T^{2q}}{q^2}N^2_2
	\left(\int_0^te^{-2\omega(t-s)\alpha}ds\right)^\frac{1}{\alpha}
	\left(\int_0^tL_\sigma^{\beta}(s)ds\right)^\frac{1}{\beta}\mathbb{E}\|x(t)-y(t)\|^2\\
	&+3\frac{T^{2q-1}}{2q-1}N^2_2\overline{c_p}\left(\int_0^te^{-2\omega(t-s)\alpha}ds\right)^\frac{1}{\alpha}\left(\int_0^tL_g^{\beta}(s) ds\right)^\frac{1}{\beta}\mathbb{E}\|x(t)-y(t)\|^2\\
	\leq& 3\frac{T^{2q-1}}{2q-1}N_2^2\left[\frac{1-e^{-2\omega T\alpha }}{2\omega\alpha}\right]^\frac{1}{\alpha}\|L_f\|_{L^{\beta}(J,\mathbb{R}^+)}\mathbb{E}\|x(t)-y(t)\|^2\\
	&+3\frac{T^{2q}}{q^2}N_2^2
	\left[\frac{1-e^{-2\omega T\alpha }}{2\omega\alpha}\right]^\frac{1}{\alpha}\|L_\sigma\|_{L^{\beta}(J, \mathbb{R}^+)}\mathbb{E}\|x(t)-y(t)\|^2
	\end{align*}\begin{align*}
	&+3\frac{T^{2q-1}}{2q-1}N_2^2\overline{c_p}\left[\frac{1-e^{-2\omega T\alpha}}{2\omega\alpha}\right]^\frac{1}{\alpha}\|L_g\|_{L^{\beta}(J,\mathbb{R}^+)}\mathbb{E}\|x(t)-y(t)\|^2\\
	\leq& 3N_2^2\left[\frac{1-e^{-2\omega T\alpha }}{2\omega\alpha}\right]^\frac{1}{\alpha}\left(\frac{T^{2q-1}}{2q-1}\|L_f\|_{L^{\beta}(J,\mathbb{R}^+)}+\frac{T^{2q}}{q^2}\|L_\sigma\|_{L^{\beta}(J,\mathbb{R}^+)}+\frac{T^{2q-1}}{2q-1}\overline{c_p}\|L_g\|_{L^{\beta}(J,\mathbb{R}^+)}\right)\\&
	\quad\times\mathbb{E}\|x(t)-y(t)\|^2 \\
	\leq& M \mathbb{E}\|x(t)-y(t)\|^2.
	\end{align*}
	Thus, $\Phi$ is a contraction mapping, it has a unique fixed point $x(t)\in\mathcal{K}_r,$ which is a solution of (\ref{p}).

	Now, we prove the stability  results of system (\ref{p}).
	For any given $\epsilon>0$, there exist $$\delta=\frac{\epsilon(1-Q_1)-Q_2}{4N^2_1e^{-2\omega T}}$$ such that $\mathbb{E}\|x_0\|\leq\delta$ implies that
	\begin{align*}
	\mathbb{E}\|x(t)\|^2
	\leq& 4N_1^2 e^{-2\omega T}\mathbb{E}\|x_0\|^2+8N_2^2r\left[\frac{1-e^{-2\omega T\alpha }}{2\omega\alpha}\right]^\frac{1}{\alpha}\\
	&\times
	\Bigg(\frac{T^{2q-1}}{2q-1}\|L_f\|_{L^{\beta}(J,\mathbb{R}^+)}
	+\frac{T^{2q}}{q^2}\|L_\sigma\|_{L^{\beta}(J,\mathbb{R}^+)}+\frac{T^{2q-1}}{2q-1}\overline{c_p}\|L_g\|_{L^{\beta}(J,\mathbb{R}^+)}\Bigg)\\
	&+8N_2^2\left[\frac{1-e^{-2\omega T}}{2\omega}\right]\Bigg(\frac{T^{2q-1}}{2q-1}R_f
	+\frac{T^{2q}}{q^2}R_\sigma+\frac{T^{2q-1}}{2q-1}\overline{c_p}R_g\Bigg)\\
	r\leq& 4N_1^2 e^{-2\omega T}\mathbb{E}\|x_0\|^2+Q_1r+Q_2\\
	r-Q_1r\leq& 4N^2_1e^{-2\omega T}\|x_0\|^2+Q_2\\
	r(1-Q_1)\leq& 4N^2_1e^{-2\omega T}\delta+Q_2\\
	r\leq&\frac {4N^2_1e^{-2\omega T}\delta+Q_2}{(1-Q_1)}\\
	r\leq&\epsilon.
	\end{align*}
	This complete the proof.
\end{theorem}
\begin{theorem}
	Assume that hypotheses (H2) - (H3) hold. Then, the system (\ref{p}) is exponential stable provided that
	\begin{align*}
	\omega>\beta=N_2^2
	\Bigg[\frac{T^{2q-1}}{2q-1}\Big(\Big[\widehat{V}_f+\overline{c_p}\widehat{V}_g\Big](1+r)\Big)+\frac{T^{2q}}{q^2}\widehat{V}_\sigma(1+r)\Bigg].
	\end{align*}
	\textbf{Proof:}
	\begin{align*}
	\mathbb{E}\|x(t)\|^2\leq&4N_1^2e^{-2\omega t}\mathbb{E}\|x_0\|^2+ 4N_2^2
	\Bigg[\frac{T^{2q -1}}{2q-1}(1+r)\Big[\widehat{V}_f+\overline{c_p}\widehat{V}_g\Big]+\frac{T^{2q}}{q^2}\widehat{V}_\sigma(1+r)\Bigg] e^{-2\omega t}\\
	&\times\int_{0}^te^{2\omega s}ds\\
	e^{2\omega t}\mathbb{E}\|x(t)\|^2
	&\leq 4N_1^2 \mathbb{E} \|x_0\|^2+4N_2^2
	\Bigg[\frac{T^{2q-1}}{2q-1}(1+r)\Big[\widehat{V}_f+\overline{c_p}\widehat{V}_g\Big]+\frac{T^{2q}}{q^2}\widehat{V}_\sigma(1+r)\Bigg]\int_0^te^{2\omega s}ds.
	\end{align*}
	Using Gronwall's inequality, we obtain
	\begin{align*}
	&e^{2\omega t}\mathbb{E}\|x(t)\|^2\leq 4N_1^2 \mathbb{E} \|x_0\|^2e^{\Bigg(4N_2^2
		\Bigg[\frac{T^{2q-1}}{2q-1}\Big(\Big[\widehat{V}_f+\overline{c_p}\widehat{V}_g\Big](1+r)\Big)+\frac{T^{2q}}{q^2}\widehat{V}_\sigma(1+r)\Bigg]t\Bigg)}.
	\end{align*}
	Consequently,
	\begin{align*}
	\mathbb{E}\|x(t)\|^2\leq& M_1 \mathbb{E}\|x_0\|^2 e^{(-\frac{v}{2}t)}.
	\end{align*}
	where $\frac{v}{2}=\omega-2\beta,$ $M_1=4N_1^2.$
	This complete the proof.
\end{theorem}
\section{Example} Set $T=1, q=0.6$ and step size $h=0.01.$
Consider the following equations
\begin{align}
^CD_{0+}^q x_1(t)=& -0.1x_1(t)-\frac{ x_2^2(t)}{1-t}-\sigma_1 x_1(t)dB_1-\frac{(0.2-t)x_1(t)e^{-t}}{\eta}\label{33}
\\\label{44}
^CD_{0+}^q x_2(t)=& -0.1x_2(t)-\frac{x_1^2(t)}{1-t}-\sigma_2 x_2(t)dB_2-\frac{(0.2-t)x_2(t)e^{-t}}{\eta},\ t\in J_1=[0,1] ,	
\end{align}
where $$ A=
\begin{pmatrix}
-0.1& 0\\
0& -0.1
\end{pmatrix},~~~
f(t,x(t))=
\begin{pmatrix}
\frac{-x_2^2(t)}{1-t}\\
\frac{-x_1^2(t)}{1-t}
\end{pmatrix},~~~
\sigma(t,x(t))dw(t)=
\begin{pmatrix}
-\sigma_1 x_1(t)d\omega_1\\
-\sigma_2 x_2(t)d\omega_2\\
\end{pmatrix},$$
$$g(t,x(t),\eta)=
\begin{pmatrix}
\frac{-(0.2-t)x_1(t)e^{-t}}{\eta}\\
\frac{-(0.2-t)x_2(t)e^{-t}}{\eta}
\end{pmatrix},~~~
\sigma_1=9.8,~~~\sigma_2=10,~~~\eta=1.$$
Now, we need to verify the hypothesis (H1) to estimate equation (\ref{22}). Let us take for any $x(t),y(t)\in \mathbb{R}^2$ and $t\in J_1$ and by simple calculation, one can obtain the following inequalities 
\begin{itemize}
	\item $\mathbb{E}\|f(t,x)-f(t,y)\|^2
	\leq -t ~\mathbb{E}\| x- y\|^2$
	\item $\mathbb{E}\|\sigma(t,x)-\sigma(t,y)\|^2
	\leq -10 t~\mathbb{E}\| x- y\|^2$
	\item  $\int_{-\infty}^{+\infty}\mathbb{E}\|g(t,x,\eta)-g(t,y,\eta)\|^2\varPi(d\eta)
	\leq -(0.2-t)~\mathbb{E}\| x-y\|^2.$
\end{itemize}
Hence, $f,\sigma$ and $g$ satisfy  the hypothesis (H1), where we set $L_f(\cdot)=-(\cdot)\in L^{\beta}(J_1,\mathbb{R}^+),L_\sigma(\cdot)=-10(\cdot)\in L^{\beta}(J_1,\mathbb{R}^+)$ and $L_g(\cdot)=-(0.2-\cdot)\in L^{\beta}(J_1,\mathbb{R}^+).$ By simple calculation, one can get  
$\|L_f\|_{L^2{(J_1,\mathbb{R}^+)}}=-1,$
$\|L_\sigma\|_{L^2{(J_1,\mathbb{R}^+)}}= -10,$
$\|L_g\|_{L^2{(J_1,\mathbb{R}^+)}}=-0.40988,$ $\|A\|=0.1000$ and $N_2=1.0202.$ Next, when we choose $\alpha=\beta=2,$ we get 
\begin{align*}
M=&\left\{3\left(\frac{1^{2\times 0.6}}{0.5^2}\right)(1.0202)(1.10885)(-11.4098)\right\}\\
=& (12.2424)(1.10885)(-11.4098)\\
=&-16.5630<1.
\end{align*}
Which guarantees that the equation (\ref{22}) holds. Thus, all the hypotheses of Theorem \ref{th1} are satisfied. Hence the equations  (\ref{33})-(\ref{44}) are stable on $J_1.$ The corresponding stability results of the system (\ref{33})-(\ref{44}) are depicted in the following figures for various fractional order $q=0.6,~ 0.75$ and $0.9.$  Further, all hypotheses of nonlinear functions $f,\sigma$ and $g$ are verified numerically. Hence, by Theorem \ref{th1}, the system (\ref{33})-(\ref{44}) is stable on $[0,100].$
\begin{figure}[h!]
	\centering
	\includegraphics[width=1.1\linewidth]{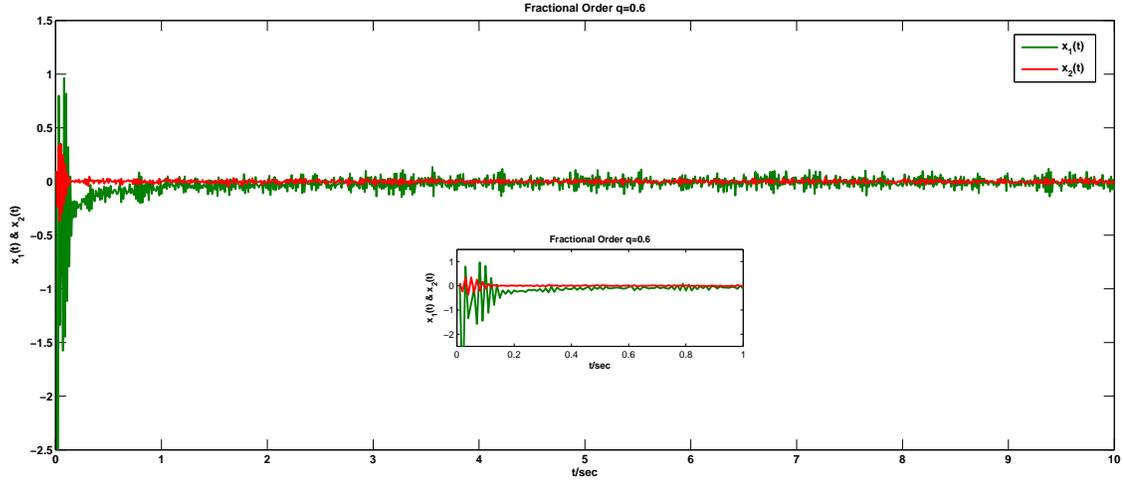}
	\caption{Time response of the states $x_1(t)\ \&\ x_2(t)$ of the system (\ref{33})-(\ref{44}) for  $q=0.6$}
	\label{fig:33}
\end{figure}
\newpage
\begin{figure}[h!]
	\centering
	\includegraphics[width=1.1\linewidth]{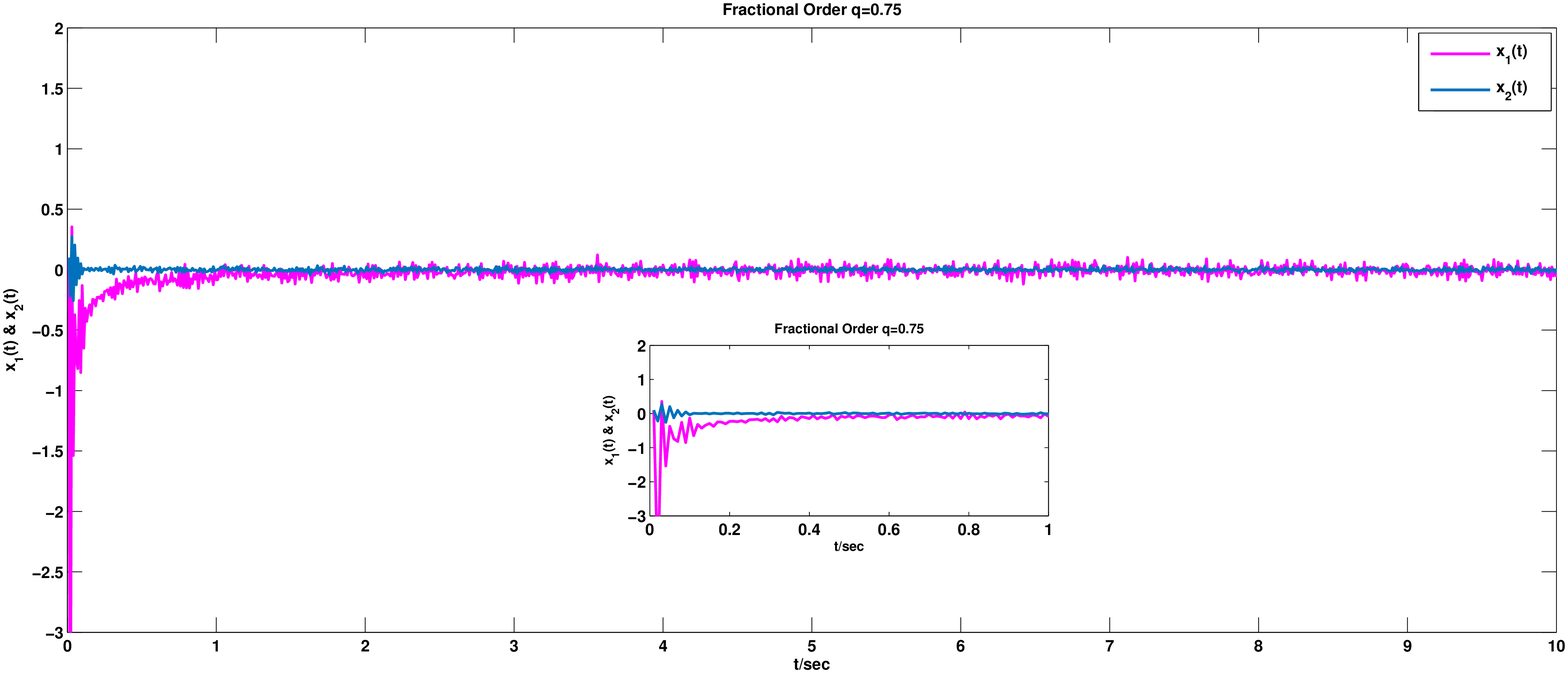}
	\caption{Time response of the states $x_1(t)\ \&\ x_2(t)$ of the system (\ref{33})-(\ref{44}) for $q=0.75$}
	\label{fig:5}
\end{figure}
\begin{figure}[h!]
	\centering
	\includegraphics[width=1.1\linewidth]{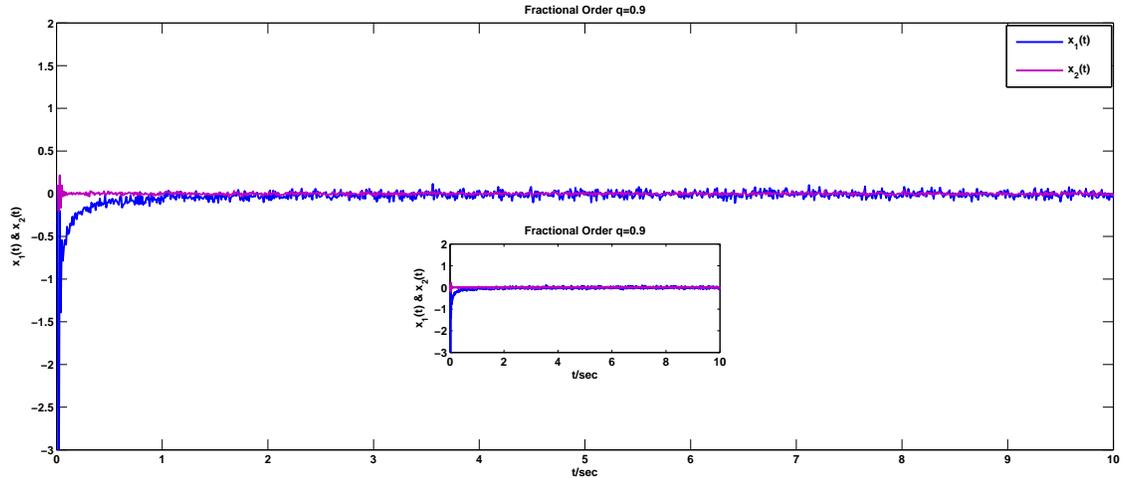}
	\caption{Time response of the states $x_1(t)\ \&\ x_2(t)$ of the system (\ref{33})-(\ref{44}) for $q=0.9$}
	\label{fig:7}
\end{figure}
\begin{remark}
	The obtained theoretical results are verified through numerical simulation. Here, state response of the proposed equations is depicted in Figures 1-3 for various fractional order $q=0.6,0.75,0.9.$ The small figures at in side of the Figures 1-3 represent the stochastic nature of the time response of the states. Hence, we conclude that the proposed system is exponentially stable.	
\end{remark}	
\section{Conclusion}
Exponential stability of nonlinear fractional order stochastic system with Poisson jumps is promoted. Main results have been obtained based on bounded properties of Mittag-Leffler function, Banach contraction mapping principle and stability theory. Finally, a numerical example have been provided to show the effectiveness of the obtained results.

\bibliographystyle{plain}



\end{document}